\documentclass[12pt]{amsart}
\usepackage{graphicx, overpic, epsf,amssymb,amscd,amsfonts}

\newtheorem{thm}{Theorem}

\newtheorem{lemma}[thm]{Lemma}

\newcommand{\R}{\mathbb{R}}
\newcommand{\Z}{\mathbb{Z}}
\newcommand{\Q}{\mathbb{Q}}

\newcommand{\bdry}{\partial}

\newcommand{\be}{\begin{enumerate}}
\newcommand{\ee}{\end{enumerate}}

\begin{document}

\title{The vanishing of the contact invariant in the presence of
torsion}

\author{Paolo Ghiggini}
\address{California Institute of Technology, Pasadena, CA 91125}
\email{ghiggini@caltech.edu}
\urladdr{http://www.its.caltech.edu/~ghiggini}


\author{Ko Honda}
\address{University of Southern California, Los Angeles, CA 90089} \email{khonda@usc.edu}
\urladdr{http://rcf.usc.edu/\char126 khonda}

\author{Jeremy Van Horn-Morris}
\address{University of Texas, Austin, TX 78712}
\email{jvanhorn@math.utexas.edu}
\urladdr{http://www.ma.utexas.edu/\char126 jvanhorn}

\date{This version: January 16, 2008.}



\thanks{This work was started during the workshop ``Interactions of
Geometry and Topology in Low Dimensions'' at BIRS. PG was supported
by a CIRGET fellowship and by the Chaire de Recherche du Canada en
alg\`ebre, combinatoire et informatique math\'ematique de l'UQAM. KH
was supported by an NSF CAREER Award (DMS-0237386).}

\begin{abstract}
We prove that the Ozsv\'ath-Szab\'o contact invariant of a closed
3-manifold with positive $2\pi$--torsion vanishes.
\end{abstract}

\maketitle

In 2002, Ozsv\'ath and Szab\'o~\cite{OSz} defined an invariant of a
closed, oriented, contact 3-manifold $(M,\xi)$ as an element of the
Heegaard Floer homology group $\widehat{HF}(-M)$. The definition of
the contact invariant was made possible by the work of
Giroux~\cite{Gi3}, which related contact structures and open book
decompositions.  The Ozsv\'ath-Szab\'o contact invariant has
undergone an extensive study, e.g., \cite{LS1,LS2}. Recently, Honda,
Kazez and Mati\'c~\cite{HKM3} defined an invariant of a contact
3-manifold with convex boundary as an element of Juh\'asz' sutured
Floer homology~\cite{Ju1,Ju2}.  The goal of this paper is to use
this relative contact invariant to prove a vanishing theorem in the
presence of torsion.

Recall that a contact manifold $(M,\xi)$ has {\em positive $n\pi
$-torsion} if it admits an embedding $(T^2\times[0,1],\eta_{n\pi
})\hookrightarrow (M,\xi)$, where $(x,y,t)$ are coordinates on
$T^2\times [0,1]\simeq \R^2/\Z^2\times[0,1]$ and $\eta_{n\pi}=\ker
(\cos (n\pi t) dx - \sin (n\pi t) dy)$.  The torsion was an
essential ingredient for distinguishing tight contact structures on
toroidal 3-manifolds (see for example~\cite{Gi1}), and is a source
of non-finiteness of the number of isotopy classes of tight contact
structures (\cite{CGH,Co,HKM1}).

\begin{thm}[Vanishing Theorem]\label{torsion}
If a closed contact 3-manifold $(M,\xi)$ has positive $2
\pi$-torsion, then its contact invariant $c(M,\xi)$ in
$\widehat{HF}(-M)$ vanishes.
\end{thm}

The coefficient ring of $\widehat{HF}(-M)$ is $\Z$ in
Theorem~\ref{torsion}. The behavior of the contact invariant with
twisted coefficients in presence of torsion is the subject of a
forthcoming paper by the first two authors \cite{GH}.

Theorem~\ref{torsion} was first conjectured in \cite[Conjecture
8.3]{Gh2}, and partial results were obtained by \cite{Gh1},
\cite{Gh2}, and \cite{LS3}. The corresponding vanishing result for
the contact class in monopole Floer homology has recently been
announced by Mrowka and Rollin (and is motivated by \cite{Ga}).
Theorem \ref{torsion}, together with a non-vanishing result of the
contact invariant proved by Ozsv\'ath and Szab\'o \cite[Theorem
4.2]{OSz:genus}, implies that a contact manifold with positive $2
\pi$-torsion is not strongly symplectically fillable. This
non-fillability result was conjectured by Eliashberg, and first
proved by Gay \cite{Ga}.

In this paper, a contact structure $\xi$ on a compact, oriented
$3$-manifold $N$ with convex boundary $\bdry N$ and dividing set
$\Gamma$ on $\bdry N$ will be denoted $(N,\Gamma,\xi)$. We will
write the invariant for a closed contact $3$-manifold $(M,\xi)$ as
$c(M,\xi)\in \widehat{HF}(-M)$ and the invariant for a compact
contact 3-manifold $(N,\Gamma,\xi)$ as $c(N,\Gamma,\xi)\in
SFH(-N,-\Gamma)$, where $SFH(-N,-\Gamma)$ is the sutured Floer
homology of $(-N,-\Gamma)$, and $\Gamma\subset \bdry N$ is now
viewed as a balanced suture. Strictly speaking, the contact
invariants have a $\pm 1$ ambiguity, but this will not complicate
matters in this paper. The key property of the relative contact
invariant which we use in this paper is the following theorem
from~\cite{HKM3}:

\begin{thm} [{\cite[Theorem 4.5]{HKM3}}] \label{thm:restriction}
Let $(M,\xi)$ be a closed contact 3-manifold and $N\subset M$ be a
compact submanifold $($without any closed components$)$ with convex
boundary and dividing set $\Gamma$. If $c(N,\Gamma, \xi|_N)=0$, then
$c(M,\xi)=0$.
\end{thm}

The behavior of the contact invariant under contact surgery will
also play a fundamental role in the proof of Theorem \ref{torsion}.

\begin{lemma}\label{surgery}
If $(N',\Gamma', \xi')$ is obtained by contact $(+1)$-surgery on a
Legendrian curve in $(N,\Gamma, \xi)$, then the contact
$(+1)$-surgery gives rise to a natural map:
\begin{equation}\label{equation: natural}
\Phi \colon  SFH(-N,-\Gamma)  \to SFH(-N',-\Gamma'),
\end{equation}
which satisfies $\Phi(c(N, \Gamma, \xi))= c(N', \Gamma', \xi')$.
\end{lemma}

\begin{proof}
If $(N', \xi')$ is obtained from $(N, \xi)$ by contact
$(+1)$-surgery, then $(N, \xi)$  is obtained from $(N', \xi')$ by
contact $(-1)$-surgery (i.e., Legendrian surgery); see
\cite[Proposition 8]{DG1}. The proof that the contact invariant is
natural with respect to Legendrian surgery is the same as in the
closed case, provided we use the reformulation of the contact
invariant given by Honda, Kazez, and Mati\'c~\cite{HKM2}.  The proof
in the closed case is given in \cite[Proposition 3.7]{HKM2}. See
also \cite[Proposition 4.4]{HKM3}.
\end{proof}

In this paper we assume that the reader is familiar with the
terminology introduced in \cite{H1}, such as {\em basic slice}, {\em
standard neighborhood of a Legendrian curve}, {\em Legendrian ruling
curve}, and {\em minimally twisting}.

Let $\Gamma$ be the following suture/dividing set on the boundary of
$T^2\times[0,1]$: $\#\Gamma_{T_0}=\#\Gamma_{T_1}=2$,
$\mbox{slope}(\Gamma_{T_0})= -1$, and
$\mbox{slope}(\Gamma_{T_1})=-2$. Here $\#$ denotes the number of
connected components, $T_i=T^2\times\{i\}$, the slope is calculated
with respect to a fixed oriented identification $T^2\simeq
\R^2/\Z^2$, and the orientation of $T_i$ is inherited from that of
$T^2$. (Hence $\bdry (T^2\times[0,1])=T_1\cup -T_0$.)

Let $\zeta_0$ be a tight contact structure so that
$(T^2\times[0,1],\Gamma,\zeta_0)$ is a basic slice. There are two
possible isotopy classes rel boundary, and $\zeta_0$ can be in
either one.

\begin{lemma} \label{embedding}
Let $L$ be a Legendrian ruling curve with infinite slope on a
parallel copy $T_\varepsilon$ of $T_0$ with the same dividing set,
inside the basic slice $(T^2\times[0,1],\Gamma, \zeta_0)$. Then
there is an embedding $i$ of $(T^2 \times [0,1],\Gamma,\zeta_0)$
into the standard tight $(S^3,\xi_{std})$, so that $i(L)$ is an
unknot with Thurston-Bennequin invariant $-1$.
\end{lemma}

\begin{figure}\centering
\includegraphics[width=6cm]{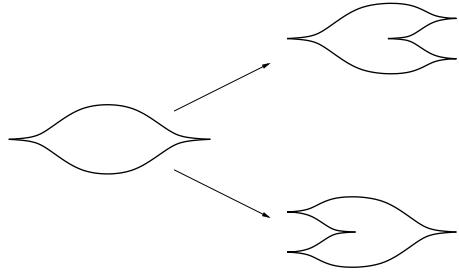}
\caption{Positive and negative stabilizations of the Legendrian unknot in $S^3$
with Thurston--Bennequin number $-1$.}
\label{stabilization}
\end{figure}

\begin{proof}
Choose coordinates $(x,y)$ on $T^2\simeq \R^2/\Z^2$ and $z$ on
$[0,1]$.  Then $(T^2 \times [0,1],\Gamma,\zeta_0)$ is contact
isomorphic to a basic slice with boundary slopes $- \frac 12$ and
$-1$ via the diffeomorphism $(x,y,z) \mapsto (y,x,1-z)$. Under this
diffeomorphism $L$ is mapped to a curve with slope $0$.

Let $V$ be a standard neighborhood of a Legendrian unknot $K$ in
$(S^3,\xi_{std})$ with Thurston-Bennequin number $-1$.  Then
$\mbox{slope}(\Gamma_{\bdry V})=-1$ and $\#\Gamma_{\bdry V}=2$,
where the slopes are computed with respect to a basis on
$H_1(\partial V)$ such that the meridian has slope $0$ and the
longitude determined by the Seifert surface has slope $\infty$. If
we stabilize $K$ and let $V'$ be a sufficiently small standard
neighborhood of the stabilized curve, then the collar region $V
\setminus V'$ is a basic slice with boundary slopes $- \frac 12$ and
$-1$. Recall that $K$ can be stabilized in two different ways, which
correspond to two different basic slices --- it is easy to relate
the relative Euler class of the basic slice with the rotation number
of the stabilized knot. See Figure \ref{stabilization} for the two
different stabilizations of $K$, drawn in the front projection.

The basic slice $(T^2 \times [0,1],\Gamma, \zeta_0)$ with boundary
slopes $-1$ and $-2$ is contact isomorphic to the basic slice $(V
\setminus V', \Gamma_{\bdry V'} \cup \Gamma_{\bdry V}, \xi_{std}|_{V
\setminus V'})$, and the Legendrian knot $L \subset T^2\times[0,1]$
corresponds to a pushoff of the meridian of $V$. Therefore, the
image of $L$ is an unknot, and the Thurston-Bennequin invariant is
easily calculated from the number of intersections with
$\Gamma_{\bdry V}$.
\end{proof}

\begin{lemma} \label{little}
The contact manifold $(M,\xi)$ has positive $2 \pi$-torsion if and
only if there exists an embedding of $(T^2\times[0,1],\Gamma,
\zeta_1)$ into $(M,\xi)$, where $(T^2\times[0,1],\Gamma, \zeta_1)$
is not minimally twisting and is homotopic relative to the boundary
to a basic slice $(T^2\times[0,1],\Gamma, \zeta_0)$.
\end{lemma}

\begin{proof}
From the classification of tight contact structures on
$T^2\times[0,1]$ (see Theorem~2.2 as well as the discussion in
Section~5.2 in \cite{H1}; an equivalent result is given in
\cite{Gi2}) it follows that, if $\zeta_1$ is not minimally twisting
and is homotopic to a basic slice, then $(T^2\times[0,1], \zeta_1)$
has positive $2 \pi$-torsion. Therefore the existence of an
embedding of $(T^2\times[0,1], \zeta_1)$ into $(M,\xi)$ implies that
$(M,\xi)$ has positive $2 \pi$-torsion.

Assume $(M, \xi)$ contains a contact submanifold isomorphic to
$(T^2\times[0,1], \eta_{2\pi})$. Then it also contains a slightly
larger submanifold $(N, \zeta')$, where $N= T^2\times[-
\varepsilon_0 ,1+ \varepsilon_1 ]$, and $\zeta'$ is defined by the
same contact form as $\eta_{2\pi}$. This can be easily seen from the
normal form of a contact structure in the neighborhood of a
pre-Lagrangian torus. By direct computation, we can choose
$\varepsilon_0, \varepsilon_1\geq 0$ so that the tori $T^2\times \{
- \varepsilon_0 \}$ and $T^2\times \{ 1+ \varepsilon_1 \}$ are
pre-Lagrangian tori with rational slopes $s_1$, $s_2$ forming an
integer basis of $H_1(T^2)$. Then we can perturb the boundary of $N$
to make it convex, so that the boundary tori have $\#\Gamma=2$ and
slopes $s_1$, $s_2$; see for example \cite[Lemma~3.4]{Gh3}. Let
$\zeta_1$ be the resulting contact structure: the contact manifold
$(N, \zeta_1)$ constructed in this way is clearly
non-minimally-twisting. After a change of coordinates in $N$, we can
make its boundary slopes $-1$ and $-2$. The contact structure is
homotopic to a basic slice by a standard explicit computation (see
\cite[Proposition 6.1]{Gh2}).
\end{proof}

\begin{proof}[Proof of Theorem~\ref{torsion}.]
By Theorem~\ref{thm:restriction} and Lemma~\ref{little}, it suffices
to prove that $c(N,\Gamma,\zeta_1)=0$, where $N=T^2\times[0,1]$ and
$\Gamma$, $\zeta_1$ are as defined above.  This proof is modeled on
the argument in \cite{Gh2}.

Take a parallel copy $T_{\varepsilon}$ of $T_0$ in the interior of
$N$ with the same dividing set, and let $L$ be a Legendrian ruling
curve on $T_{\varepsilon}$ with slope $\infty$. The Legendrian curve
$L$ has twisting number $-1$ with respect to the framing coming from
$T_{\varepsilon}$.

Now apply a contact $(+1)$-surgery to $N$ along $L$; see for
example \cite{DG2}. As the surgery coefficient is $0$ with respect
to the framing induced by the torus $T_{\varepsilon}$, the resulting
3-manifold is $N'=(S^1\times D^2) \# (S^1\times D^2)$. Next write
$\Gamma'$ as $\Gamma'_1\sqcup \Gamma_2'$, where $\Gamma_i'$ is the
dividing set on the $i$th connect summand $S^1\times D^2$.  Since
each component of $\Gamma_i'$ intersects the meridian once
geometrically, we may take $\Gamma_i'$ to have slope $\infty$, after
diffeomorphism. (Here the slope of the boundary of a solid torus is
defined by setting the meridian to have slope $0$ and choosing some
longitude.)

It was proved in \cite{HKM3} that $SFH(-N,-\Gamma)=\Z\oplus
\Z\oplus\Z\oplus \Z$, where each $\Z$-summand corresponds to a
distinct relative Spin$^c$-structure.  As for $SFH(-N',-\Gamma')$,
Juh\'asz~\cite[Proposition 9.15]{Ju1} proved that the sutured Floer
homology of a connected sum of two balanced sutured manifolds is the
sutured Floer homology of their tensor product, tensored with an
extra $\Z^2$ factor. Since each $(S^1\times D^2,\Gamma_i')$ is
product disk decomposable, $SFH(S^1\times D^2,\Gamma_i')\cong\Z$,
and hence $SFH(-N',-\Gamma')\cong (\Z\otimes \Z) \otimes \Z^2 \cong
\Z^2$.

Let $\mathfrak{s}$ be the relative Spin$^c$-structure induced by
$\zeta_1$. We claim that the map $\Phi$ induced by the surgery is injective
on the direct summand $SFH(-N,-\Gamma, \mathfrak{s}) \cong \Z$; that is the
content of Lemma~\ref{injectivity} below. In
Lemma~\ref{overtwisted}, we will prove that applying contact
$(+1)$-surgery to $(N,\Gamma,\zeta_1)$ along $L$ yields an
overtwisted contact structure $\zeta_1'$ on $N'$.  Therefore,
$\Phi (c(N,\Gamma,\zeta_1))= c(N',\Gamma',\zeta_1')=0$, and by the injectivity of $\Phi$ on the
appropriate $\Z$-summand it follows that $c(N,\Gamma,\zeta_1)=0$.
\end{proof}

\begin{lemma} \label{injectivity}
Let $\mathfrak{s}$ be the relative Spin$^c$-structure induced by
$(\Gamma,\zeta_1)$ and $\mathfrak{s'}$ be that induced by
$(\Gamma',\zeta_1')$. Then the map
\[ \Phi \colon SFH(-N,-\Gamma, \mathfrak{s}) \to SFH(-N',-\Gamma', \mathfrak{s'}), \]
given by Equation~\ref{equation: natural}, is injective.
\end{lemma}

\begin{proof}
Recall that $\zeta_0$ and $\zeta_1$ have the same relative
Spin$^c$-structure $\mathfrak{s}$. By
Lemma~\ref{embedding}, $(N,\Gamma,\zeta_0)$ can be embedded in
$(S^3,\xi_{std})$, which has nonzero contact invariant. Hence, by
Theorem~\ref{thm:restriction}, the contact invariant
$c(N,\Gamma,\zeta_0)\in SFH(-N,-\Gamma, \mathfrak{s})$ is nonzero.
Since $SFH(-N,-\Gamma, \mathfrak{s})\cong \Z$ (since it is nonzero)
and $SFH(-N',-\Gamma')\cong \Z^2$, it suffices to prove that
$\Phi(c(N,\Gamma,\zeta_0)) \neq 0$.

By Lemma \ref{surgery}, the cobordism map $\Phi$ takes the contact
class $c(N,\Gamma,\zeta_0)$ to $c(N',\Gamma',\zeta_0')$, where
$\zeta_0'$ is the contact structure obtained from $\zeta_0$ by
contact $(+1)$-surgery along $L$. Now consider the embedding
$i \colon (N,\Gamma,\zeta_0)\hookrightarrow (S^3,\xi_{std})$ from
Lemma~\ref{embedding}. Legendrian $(+1)$-surgery along the unknot
$i(L)$ with Thurston-Bennequin invariant $-1$ inside
$(S^3,\xi_{std})$ yields the unique tight contact structure on
$S^1\times S^2$, which has nonzero contact invariant: for example,
see \cite[Lemma 3.7]{LS2}. Hence $c(N',\Gamma',\zeta_0')\not=0$, and
it follows that $SFH(-N,-\Gamma,\mathfrak{s})$ maps injectively into
$SFH(-N',-\Gamma')$.
\end{proof}

\begin{lemma} \label{overtwisted}
Applying contact $(+1)$-surgery to $(N,\Gamma,\zeta_1)$ along $L$
yields an overtwisted contact structure $\zeta_1'$ on $N'$.
\end{lemma}

\begin{proof}
For any $s \in \Q \cup \{ \infty \}$, there is a convex torus (in
standard form) with slope $s$ in $(N,\Gamma,\zeta_1)$ parallel to
the boundary, according to \cite[Proposition 4.16]{H1}. In
particular, there is a standard torus whose Legendrian divides have
the same slope as the Legendrian ruling curve $L$ we are doing
surgery on. After the surgery, this Legendrian divide bounds an
overtwisted disk in $N'$.
\end{proof}


\begin{thebibliography}{}

\bibitem[Co]{Co}
V.\ Colin, \textit{Une infinit\'e de structures de contact tendues
sur les vari\'et\'es toro\"idales}, Comment.\ Math.\ Helv.\ {\bf 76}
(2001), 353--372.

\bibitem[CGH]{CGH}
V.\ Colin, E.\ Giroux and K.\ Honda, \textit{Finitude homotopique et
isotopique des structures de contact tendues}, in preparation.

\bibitem[DG1]{DG1}
F.\ Ding and H.\ Geiges, \textit{Symplectic fillability of tight contact
structures on torus bundles}, Alg. \ Geom. \ Topol. \ {\bf 1} (2001),
153--172.

\bibitem[DG2]{DG2}
F.\ Ding and H.\ Geiges, \textit{A Legendrian surgery presentation
of contact 3-manifolds},  Math.\ Proc.\ Cambridge Philos.\ Soc.\
{\bf 136} (2004), 583--598.

\bibitem[Ga]{Ga}
D.\ Gay, \textit{Four-dimensional symplectic cobordisms containing
three-handles},  Geom.\ Topol.\ {\bf 10}  (2006), 1749--1759
(electronic).

\bibitem[Gh1]{Gh1}
P.\ Ghiggini, \textit{Ozsv\'ath--Szab\'o invariants and fillability
of contact structures}, Math.\ Z.\ 253 (2006), 159--175.

\bibitem[Gh2]{Gh2}
P.\ Ghiggini, \textit{Infinitely many universally tight contact
manifolds with trivial Ozsv\'ath-Szab\'o contact invariants}, Geom.\
Topol.\ {\bf 10}  (2006), 335--357 (electronic).

\bibitem[Gh3]{Gh3}
P.\ Ghiggini, \textit{Tight contact structures on Seifert manifolds
over $T^2$ with one singular fibre},  Algebr.\ Geom.\ Topol.\ {\bf
5} (2005), 785--833 (electronic).

\bibitem[Gi1]{Gi1}
E.\ Giroux, \textit{Une infinit\'e de structures de contact tendues
sur une infinit\'e de vari\'et\'es}, Invent.\ Math.\ {\bf 135}
(1999), 789--802.

\bibitem[Gi2]{Gi2}
E.\ Giroux, \textit{Structures de contact en dimension trois et
bifurcations des feuilletages de surfaces}, Invent.\ Math.\ {\bf
141} (2000), 615--689.

\bibitem[Gi3]{Gi3}
E.\ Giroux, \textit{G\'eom\'etrie de contact:\ de la dimension trois
vers les dimensions sup\'erieures}, Proceedings of the International
Congress of Mathematicians, Vol.\ II (Beijing, 2002),  405--414,
Higher Ed.\ Press, Beijing, 2002.

\bibitem[GH]{GH}
P. Ghiggini and K. Honda, \textit{Giroux torsion and twisted coefficients},
in preparation.

\bibitem[H1]{H1}
K.\ Honda, \textit{On the classification of tight contact structures
I}, Geom.\ Topol.\ {\bf 4} (2000), 309--368 (electronic).

\bibitem[HKM1]{HKM1}
K.\ Honda, W.\ Kazez and G.\ Mati\'c, \textit{Convex decomposition
theory},  Int.\ Math.\ Res.\ Not.\ {\bf 2002}, 55--88.

\bibitem[HKM2]{HKM2}
K.\ Honda, W.\ Kazez and G.\ Mati\'c, \textit{On the contact class in Heegaard
Floer homology}, preprint 2007. \texttt{ArXiv:0609734}.

\bibitem[HKM3]{HKM3}
K.\ Honda, W.\ Kazez and G.\ Mati\'c, \textit{The contact invariant
in sutured Floer homology}, preprint 2007. \texttt{ArXiv:0705.2828}.

\bibitem[Ju1]{Ju1}
A.\ Juh\'asz, \textit{Holomorphic discs and sutured manifolds},
Algebr.\ Geom.\ Topol.\ {\bf 6} (2006), 1429--1457 (electronic).

\bibitem[Ju2]{Ju2}
A.\ Juh\'asz, \textit{Floer homology and surface decompositions},
preprint 2006. \texttt{ArXiv:math.GT/0609779}.

\bibitem[LS1]{LS1}
P.\ Lisca and A.\ Stipsicz, \textit{Ozsv\'ath-Szab\'o invariants and
tight contact three-manifolds. I}, Geom.\ Topol.\ {\bf 8} (2004),
925--945 (electronic).

\bibitem[LS2]{LS2}
P.\ Lisca and A.\ Stipsicz, \textit{Ozsv\'ath-Szab\'o invariants and
tight contact three-manifolds, II}, preprint 2004.
\texttt{ArXiv:math.SG/0404136}.

\bibitem[LS3]{LS3}
P.\ Lisca and A.\ Stipsicz, \textit{Contact Ozsv\'ath-Szab\'o
invariants and Giroux torsion},  Algebr.\ Geom.\ Topol.\ {\bf 7}
(2007), 1275--1296 (electronic).

\bibitem[OSz1]{OSz}
P.\ Ozsv\'ath and Z.\ Szab\'o, \textit{Heegaard Floer homology and
contact structures}, Duke Math.\ J.\ {\bf 129} (2005), 39--61.

\bibitem[OSz2]{OSz:genus}
P.\ Ozsv\'ath and Z.\ Szab\'o, \textit{Holomorphic disks and genus
bounds}, Geom.\ Topol.\  {\bf 8}  (2004), 311--334 (electronic).

\end{thebibliography}
\end{document}